\documentclass[10pt]{article}
\usepackage{amstext}
\usepackage{amsfonts}
\usepackage{amssymb}
\usepackage{amsbsy}
\usepackage{latexsym}
\usepackage{xy}
\usepackage{hhline}
\xyoption{all}
\newcounter{num}[section] %

\newenvironment{theo}
{\refstepcounter{num}%
\bigskip\noindent{\bf Theorem~\arabic{section}.\arabic{num}. }\it}

\newenvironment{cor}
{\refstepcounter{num}%
\bigskip\noindent{\bf Corollary~\arabic{section}.\arabic{num}. }\it}

\newenvironment{lemma}
{\refstepcounter{num}%
\bigskip\noindent{\bf Lemma~\arabic{section}.\arabic{num}. }\it}

\newcommand{\remark}
{\refstepcounter{num}%
\bigskip\noindent{\bf Remark~\arabic{section}.\arabic{num}.}}

\newcommand{\conj}
{\refstepcounter{num}%
\bigskip\noindent{\bf Conjecture~\arabic{section}.\arabic{num}.}}

\newcommand{\definition}[1]
{\refstepcounter{num}%
\bigskip\noindent{\bf Definition~\arabic{section}.\arabic{num}}~({\it #1}).}

\newenvironment{proof}{\medskip\noindent{\it Proof. }}
{$\Box$ \bigskip}
\newenvironment{proof_without_dot}{\medskip\noindent{\it Proof }}
{$\Box$ \bigskip}

\newenvironment{eq}{\begin{equation}}{\end{equation}}

\newcommand{\Ref}[1]{(\ref{#1})}

\newcommand{\si}{\sigma}
\newcommand{\al}{\alpha}
\newcommand{\be}{\beta}
\newcommand{\ga}{\gamma}
\newcommand{\la}{\lambda}
\newcommand{\de}{\delta}
\newcommand{\De}{\Delta}

\newcommand{\ov}[1]{\overline{#1}}
\newcommand{\un}[1]{{\underline{#1}} }
\newcommand{\id}[1]{{{\rm id}\{{#1}\}}}

\newcommand{\tr}{\mathop{\rm tr}}

\newcommand{\mdeg}{\mathop{\rm mdeg}}

\newcommand{\Char}{\mathop{\rm char}}

\newcommand{\algA}{\mathcal{A}}    

\newcommand{\M}{\mathcal{M}} 
\newcommand{\N}{\mathcal{N}} 
\newcommand{\bbar}[1]{\overline{#1}} 

\newcommand{\FF}{{\mathbb{F}}}   
\newcommand{\CC}{{\mathbb{C}}}   
\newcommand{\NN}{{\mathbb{N}}}

\newcommand{\DP}{{\rm DP} }

\newcommand{\Sp}{S\!p}

\begin{document}
\title{On minimal generating systems for matrix $O(3)$-invariants}
 \author{
A.A. Lopatin \\
{\small\it $Fakult\ddot{a}t$ $f\ddot{u}r$ Mathematik, %
$Universit\ddot{a}t$ Bielefeld, %
Postfach 100131, 33501 Bielefeld, Germany}  \\
{\small\it Institute of Mathematics, SBRAS, Pevtsova street, 13, Omsk 644099, Russia} \\
{\small\it artem\underline{ }lopatin@yahoo.com} \\
}
\date{} 
\maketitle

\begin{abstract} 
The algebra of invariants of several $3\times 3$ matrices under the action of the orthogonal group by simultaneous conjugation is considered over an infinite field of characteristic different from two.  The maximal degree of elements of a minimal system of generators is described with deviation $3$. 
\end{abstract}



\section{Introduction}\label{section_intro}

All vector spaces, algebras, and modules are over an infinite field $\FF$ of characteristic $\Char{\FF}\neq2$.  By an algebra we always mean an associative algebra.


The algebra of {\it matrix $O(n)$-invariants} $R^{O(n)}$ is the subalgebra of the polynomial algebra 
$$R=R_{n}=\FF[x_{ij}(k)\,|\,1\leq i,j\leq n,\, 1\leq k\leq d],$$ 
generated by $\si_t(A)$, where $1\leq t\leq n$ and $A$ ranges over all monomials in the 
matrices $X_1,\ldots,X_d$, $X_1^T,\ldots,X_d^T$. Here $\si_t(A)$ stands for the $t^{\rm th}$ coefficient of the characteristic polynomial of $A$ and 
$$X_k=\left(\begin{array}{ccc}
x_{11}(k) & \cdots & x_{1n}(k)\\
\vdots & & \vdots \\
x_{n1}(k) & \cdots & x_{nn}(k)\\
\end{array}
\right)$$
is an $n\times n$ {\it generic} matrix, where $1\leq k\leq d$. Moreover, we can assume that $A$ is {\it primitive}, i.e., is not equal to the power of a shorter monomial. By the Hilbert--Nagata Theorem on invariants, $R^{O(n)}$ is a finitely generated
algebra, but the mentioned generating system is not finite.

Similarly to matrix $O(n)$-invariants we can define matrix $GL(n)$, $\Sp(n)$ and $SO(n)$-invariants. Their generators were found in~\cite{Sibirskii_1968}, \cite{Procesi76}, \cite{Aslaksen95}, \cite{Donkin92a}, \cite{Zubkov99}, and~\cite{Lopatin_so_inv}. For $\Char{\FF}=0$ relations  between generators for $GL(n)$, $O(n)$, and $\Sp(n)$-invariants were computed in~\cite{Razmyslov74} and~\cite{Procesi76}.  For any $\Char{\FF}$ relations for $GL(n)$ and $O(n)$-invariants were established in~\cite{Zubkov96},~\cite{Lopatin_Orel}.  Note that in the case of $O(n)$-invariants we always assume $\Char{\FF}\neq2$. 

For $f\in R$ denote by $\deg{f}$ its {\it degree} and by $\mdeg{f}$ its {\it multidegree}, i.e.,
$\mdeg{f}=(t_1,\ldots,t_d)$, where $t_k$ is the total degree of the polynomial $f$ in $x_{ij}(k)$, $1\leq i,j\leq n$, and $\deg{f}=t_1+\cdots+t_d$. Since $\deg{\si_t(Y_1\cdots Y_s)}=ts$, where $Y_k$ is a generic or a transpose generic matrix, the algebra $R^{O(n)}$ has $\NN$-grading by degrees and $\NN^d$-grading by multidegrees, where $\NN$ stands for non-negative integers. 

Given an $\NN$-graded algebra $\algA$, denote by $\algA^{+}$ the subalgebra generated by homogeneous elements of positive degree. A set $\{a_i\} \subseteq \algA$ is a minimal (by inclusion) homogeneous system of generators
(m.h.s.g.)~if and only if the $a_i$'s are $\NN$-homogeneous and $\{\ov{a_i}\}$ is a basis of $\ov{\algA}={\algA}/{(\algA^{+})^2}$. If we consider $a\in\algA$ as an element of $\ov{\algA}$, then we usually omit the bar and write $a\in\ov{\algA}$ instead of $\ov{a}$. An element $a\in \algA$ is called {\it
decomposable} if $a=0$ in $\ov{\algA}$. In other words, a decomposable
element is equal to a polynomial in elements of strictly lower degree. Therefore  the highest degree of indecomposable invariants $D_{\rm max}=D_{\rm max}(n,d)$ is equal to
the least upper bound for the degrees of elements of a m.h.s.g.~for $R^{O(n)}$. So using $D_{\rm max}$ we can easily construct a finite system of generators. In this paper we give the following estimations on $D_{\rm max}$ in case $n=3$.

\begin{theo}\label{theo_appl} Let $n=3$ and $d\geq1$. Then
\begin{enumerate} 
\item[$\bullet$] If $\Char\FF =3$, then $2d+4\leq D_{\rm max}\leq 2d+7$.

\item[$\bullet$] If $\Char\FF \neq 2,3$, then $D_{\rm max}=6$.
\end{enumerate}
\end{theo}

A m.h.s.g.~for matrix $GL(2)$-invariants was found~\cite{Sibirskii_1968},~\cite{Procesi_1984}, and~\cite{DKZ_2002} and for matrix $O(2)$-invariants over $\FF=\CC$ in~\cite{Sibirskii_1968}. 
A m.h.s.g.~for matrix $GL(3)$-invariants was established in~\cite{Lopatin_Comm1} and~\cite{Lopatin_Comm2}. For $d=2$ relations for matrix $GL(3)$-invariants were explicitly described in~\cite{Nakamoto_2002},~\cite{ADS_2006}. In~\cite{Lopatin_2222} $D_{max}$ was estimated for invariants of quivers of dimension $(2,\ldots,2)$. 

\smallskip
The paper is organized as follows. In Section~\ref{section_notations} we introduce notations that are used throughout the paper. We also formulate key Theorem~\ref{theo_rel_main} from~\cite{Lopatin_Orel}.

In Section~\ref{section_A3d} we define the associative algebras $A_{3,d}$ and $N_{3,d}$ and find out some relations for them.  

In Section~\ref{section_multi} we consider relations for $\ov{R^{O(n)}}$ of multidegree $(1,\de_2\ldots,\de_d)$ for an arbitrary $n$. In Lemma~\ref{lemma1_multi} it is shown that the T-ideal of the mentioned relations is generated by two identities, whereas in Theorem~\ref{theo_rel_main} we take infinitely many identities.

In Section~\ref{section_reduction} we establish that the T-ideal of relations for $\ov{R^{O(3)}}$ is generated by four identities (see Theorem~\ref{theo_reduction} together with the definition of $A_{3,d}$). A connection between relations for $\ov{R^{O(3)}}$ and $A_{3,d}$ is described in Corollary~\ref{cor1}. Note that the algebra $A_{3,d}$ plays for $\ov{R^{O(3)}}$ a similar role as the algebra $N_{3,d}$ plays for $\ov{R^{GL(3)}}$ (see Lemmas~5 and~6 of~\cite{Lopatin_Comm1}).

In Section~\ref{section_proof} we apply results from Sections~\ref{section_A3d} and~\ref{section_reduction} to calculate $D_{max}$ and the nilpotency degree of $A_{3,d}$ with deviation 3.  Note that considering relations from Theorem~\ref{theo_reduction} of degree less or equal than $D_{\rm max}$ we obtain a {\it finite} generating system for the ideal of relations for $\ov{R^{O(3)}}$. 

Let us remark that for $n=3$ Theorem~\ref{theo_appl} implies that if $\Char\FF>n$, then $D_{\rm max}$ is the same as in the case of $\Char\FF=0$. The same phenomenon is also valid for matrix $GL(n)$-invariants for $n=2,3$. To the best of author's knowledge, there is no counterexample to the natural conjecture.

\section{Notations and known results}\label{section_notations}

For a vector $\un{t}=(t_1,\ldots,t_u)\in\NN^u$ we write $\#\un{t}=u$ and 
$|\un{t}|=t_1+\cdots +t_u$. In this paper we use the following notions from~\cite{Lopatin_Orel}: 
\begin{enumerate}
\item[$\bullet$] the monoid $\M$ (without unity) freely generated by {\it letters}  $x_1,\ldots,x_d,x_1^T,\ldots,x_d^T$, the vector space $\M_{\FF}$ with the basis $\M$, and $\N\subset\M$ the subset of primitive elements, where the notion of a primitive element is defined as above; 

\item[$\bullet$] the involution ${}^T:\M_{\FF}\to\M_{\FF}$ defined by $x^{TT}=x$ for a letter $x$ and $(a_1\cdots a_p)^T=a_p^T\cdots a_1^T$ for $a_1,\ldots,a_p\in\M$;

\item[$\bullet$] the equivalence $y_1\cdots y_p\sim z_1\cdots z_p$ that holds 
if there exists a cyclic permutation $\pi\in S_p$ such that
$y_{\pi(1)}\cdots y_{\pi(p)} = z_1\cdots z_p$ or $y_{\pi(1)}\cdots y_{\pi(p)} = z_p^T\cdots z_1^T$, where $y_1,\ldots,y_p,z_1,\ldots,z_p$ are letters;

\item[$\bullet$] $\M_{\si}$, the ring with unity of (commutative) polynomials over $\FF$ freely generated by the ``symbolic'' elements $\si_t(\al)$, where $t>0$ and $\al\in\M_{\FF}$;   

\item[$\bullet$] $\N_{\si}$, a ring with unity of (commutative) polynomials over $\FF$ freely generated by the ``symbolic'' elements $\si_t(\al)$, where $t>0$ and $\al\in\N$ ranges over $\sim$-equivalence classes; note that $\N_{\si}\simeq \M_{\si}/L$, where the ideal $L$ is described in Lemma~3.1 of~\cite{Lopatin_Orel};   

\item[$\bullet$] $\NN$-gradings by degrees and $\NN^d$-gradings by multidegrees for $\M$, $\M_{\FF}$, $\N$, and $\N_{\si}$, where $\mdeg(a)= 
(\deg_{x_1}{a}+\deg_{x_1^T}{a},\ldots,\deg_{x_d}{a}+\deg_{x_d^T}{a})$ for $a\in\M$; 

\item[$\bullet$] the surjective homomorphism $\Psi_n:\N_{\si} \to R^{O(n)}$, defined in the natural way (see Section~1 of~\cite{Lopatin_Orel}); %

\item[$\bullet$] $\si_{t,r}(a,b,c)\in\N_{\si}$ introduced in~\cite{ZubkovII}, where $a,b,c\in\M_{\FF}$ and $t,r\in\NN$ (see Section~3 of~\cite{Lopatin_Orel} or  Definition~\ref{def1} for the image of  $\si_{t,r}(a,b,c)$ in $\ov{\N_{\si}}$);

\item[$\bullet$] the partial linearization $\si_{\un{t};\un{r};\un{s}}(\un{a};\un{b};\un{c}) \in \N_{\si}$ of $\si_{t,r}(a,b,c)$ (see Section~5 of~\cite{Lopatin_Orel}), where $\un{t}=(t_1,\ldots,t_u)\in\NN^u$, $\un{r}=(r_1,\ldots,r_v)\in\NN^v$, $\un{s}=(s_1,\ldots,s_w)\in\NN^w$ satisfy $|\un{r}|=|\un{s}|$, $\un{a}=(a_1,\ldots,a_u)$, $\un{b}=(b_1,\ldots,b_v)$, $\un{c}=(c_1,\ldots,c_w)$, 
and $a_i,b_j,c_k$ belong to $\M_{\FF}$ for $1\leq i\leq u$, $1\leq j\leq v$, $1\leq k\leq w$. A formula for computation of   $\si_{\un{t};\un{r};\un{s}}(\un{a};\un{b};\un{c})$ is given in  Lemma~5.2 of~\cite{Lopatin_Orel}.
\end{enumerate} 

\noindent 
The mapping $\Psi_n$ induces the isomorphism of algebras
$$\ov{R^{O(n)}}\simeq  \ov{\N_{\si}}/\ov{K_{n}}$$
for some ideal $\ov{K_{n}}\triangleleft\ov{\N_{\si}}$. Elements of $\ov{K_{n}}$ are called {\it relations} for $\ov{R^{O(n)}}$. The ideal of relations was described in Corollary~7.6 of~\cite{Lopatin_Orel}:

\begin{theo}\label{theo_rel_main} The ideal of relations $\ov{K_{n}}$ is generated by $\si_{t,r}(a,b,c)$, where $t+2r>n$ and $a,b,c$ range over $\M_{\FF}$.
\end{theo}

\remark\label{rem_rel_main}
Since $\FF$ is infinite, we can reformulate Theorem~\ref{theo_rel_main} as follows: the ideal $\ov{K_{n}}$ is generated by the $\NN^d$-homogeneous elements 
$\si_{\un{t};\un{r};\un{s}}(\un{a};\un{b};\un{c})$, where $|\un{t}|+2|\un{r}|>n$ and 
$a_i,b_j,c_k$ range over $\M$.
\bigskip

Consider $a=y_1\cdots y_p$, where $y_1,\ldots,y_s$ are letters. We say that a letter $y_i$ is followed by $y_j$ in $a$ if $j=i+1$ or $j=1$, $i=p$. For the sake of completeness let us recall the definition of $\si_{t,r}(a,b,c)\in\ov{\N_{\si}}$.

\definition{of $\si_{t,r}(a,b,c)\in\ov{\N_{\si}}$}\label{def1} We assume $d\geq3$. Denote by $\ov{I_{t,r}}=\{a_k\}$ the set of primitive pairwise different with respect to $\sim$-equivalence words in the letters $x_i,x_i^T$, $i=1,2,3$, satisfying 
\begin{enumerate}
\item[$\bullet$] $j_{a_k}\mdeg(a_k)=(t,r,r)$ for some $j_{a_k}$;

\item[$\bullet$] the letters $x_1,x_3,x_3^T$ are followed by $x_1,x_2,x_2^T$ in $a_k$;   

\item[$\bullet$] the letters $x_1^T,x_2,x_2^T$ are followed by $x_1^T,x_3,x_3^T$ in $a_k$.  
\end{enumerate}
Then 
$$\si_{t,r}(x_1,x_2,x_3)=\sum_{a\in \ov{I_{t,r}}} (-1)^{\xi_a} \, \si_{j_a}(a) \text{ in } \ov{\N_{\si}},$$ 
where $\xi_a=t+j_a(\deg_{x_2}{a}+\deg_{x_3}{a}+1)$. We also set $\si_{0,0}(x_1,x_2,x_3)=1$. For any $d\geq1$ and $a_1,a_2,a_3\in\M_{\FF}$ we define $\si_{t,r}(a_1,a_2,a_3)\in\ov{\N_{\si}}$ as the result of the substitution $x_i\to a_i$, $x_i^T\to a_i^T$, $i=1,2,3$, in $\si_{t,r}(x_1,x_2,x_3)\in\ov{\N_{\si}}$.
\bigskip
 
The decomposition formula from~\cite{Lopatin_bplp} implies that for $n\times n$ matrices $A_i$, $i=1,2,3$, with $n=t_0+2r$,  $t_0\geq0$ we have  
\begin{eq}\label{eq1_history}
\DP_{r,r}(A_1+\la E,A_2,A_3)=\sum_{t=0}^{t_0} \la^{t_0-t}\si_{t,r}(A_1,A_2,A_3),
\end{eq}
where $\DP_{r,r}(A_1,A_2,A_3)$ stands for the determinant-pfaffian (see~\cite{LZ1}) and $\si_{t,r}(A_1,A_2,A_3)$ is defined as the result of the substitution $a_i\to A_i$, $a_i^T\to A_i^T$ in $\si_{t,r}(a_1,a_2,a_3)$. Thus $\DP_{r,r}$ relates to $\si_{t,r}$ in the same way as the determinant relates to $\si_t$.    
 
\bigskip

If $f,h\in\ov{\N_{\si}}$ are equal as elements of $\ov{R^{O(n)}}$, then we write $f\equiv h$. In particular, $f\in\ov{\N_{\si}}$ is a relation for $\ov{R^{O(n)}}$ if and only if $f\equiv0$. We say that $f\equiv h$ {\it follows} from
$f_1\equiv h_1,\ldots,f_s\equiv h_s$, if $f-h$ is a linear combination of $f_1-h_1,\ldots,f_s-h_s$ in $\ov{\N_{\si}}$. We use the following convention:
$$\tr(a)=\si_1(a)$$
for all $a\in\M_{\FF}$. Note that $\tr$ is linear in $\N_{\si}$, i.e., $\tr(\al a+\be b)=\al\tr(a)+\be\tr(b)$ in $\N_{\si}$ for $\al,\be\in\FF$ and $a,b\in\M_{\FF}$ (see Lemma~3.1 of~\cite{Lopatin_Orel}). Let us remark that 
\begin{eq}\label{eq13}
\si_{\un{t};\un{r};\un{s}}(\un{a};\un{b};\un{c}) = \si_{\un{t};\un{s};\un{r}}(\un{a}^T;\un{c};\un{b}) \text{ in }\N_{\si},
\end{eq}
where $\un{a}^T=(a_1^T,\ldots,a_u^T)$.

\section{The algebra $A_{3,d}$}\label{section_A3d}

Given $a\in\M_{\FF}$, we denote $\bar{a}=a-a^T$.  For an algebra $\algA$, denote by $\id{a_1,\ldots,a_s}$ the ideal generated by $a_1,\ldots,a_s\in\algA$. We denote

\begin{enumerate}
\item[$\bullet$] $N_{3,d}=\M_{\FF} / \id{a^3\,|\,a\in\M_{\FF}}$;

\item[$\bullet$] $T(a,b,c)=a\bar{b}\bar{c}+\bar{b}a^T\bar{c}+\bar{b}\bar{c}a= a\bar{b}\bar{c}+\bar{b}a\bar{c}+\bar{b}\bar{c}a-\bar{b}\bar{a}\bar{c}$, where $a,b,c\in\M_{\FF}$;

\item[$\bullet$] $A_{3,d}=N_{3,d}/\id{T(a,b,c)\,|\,a,b,c\in\M_{\FF}}$;

\item[$\bullet$] $\M_1=\M\sqcup\{1\}$, i.e., we endow $\M$ with the unity. 
\end{enumerate}

Since $\FF$ is infinite, the elements 
\begin{enumerate}
\item[$\bullet$] $T_1(a)=a^3$,

\item[$\bullet$] $T_2(a,b)=a^2b + aba + ba^2$,

\item[$\bullet$] $T_3(a,b,c)=abc+acb+bac+bca+cab+cba$,
\end{enumerate}
where $a,b,c\in\M_{\FF}$, are equal to zero in $N_{3,d}$. Moreover,  
$$N_{3,d} =\M_{\FF}/\id{T_1(a),\, T_2(a,b),\,T_3(a,b,c)\,|\,a,b,c\in\M}.$$

In what follows we apply the following remark without references to it.

\remark\label{rem1} Consider $f=\sum_i\al_i a_ibc_i$, where $\al_i\in\FF$, $a_i,c_i\in\M_1$, and $b\in\M$. If $f=0$ in $A_{3,d}$ is valid for all $b\in\M$, then $f=0$ in $A_{3,d}$ is also valid for all $b\in\M_{\FF}$. 
\bigskip 

The following equalities in $\M$ are trivial: %
$$\bar{\bar{a}}=2\bar{a},\; \bbar{a^{T}}=-\bar{a}=\bar{a}^T,\; 
\bbar{ab}=ab-ba+\bar{b}a+b\bar{a}-\bar{b}\bar{a}.$$ %

\begin{lemma}\label{lemma_N3}
We have
\begin{enumerate}
\item[a)] If $a\in\M$ and $x$ is a letter, then $a\in A_{3,d}$ is equal to the sum of the  following elements: $b_1$, $b_1xb_2$, $b_1x^2b_2$,
$b_1x^2 c x b_2$, where $b_1,b_2\in\M_1$, $\deg_{x}(b_i)=\deg_{x^T}(b_i)=0$ for $i=1,2$, and $c\in\M$; in particular, if $\deg_x{a}>3$, then $a=0$ in $A_{3,d}$;

\item[b)] If $\Char{\FF}=0$ or $\Char{\FF}>3$, then $x_1\cdots x_6=0$ in $A_{3,d}$;

\item[c)] $(ab)^2=b^2a^2$ in $A_{3,d}$, where $a,b\in\M$; moreover, $(a_1\cdots a_s)^2=a_s^2\cdots a_1^2$ in $A_{3,d}$, where $a_1,\ldots,a_s\in\M$;

\item[d)] $ab\cdot c\cdot ba=-ca^2 b^2 - b^2a^2c$ in $A_{3,d}$, where $a,b,c\in\M$; 

\item[e)] If $\Char{\FF}=3$ and for $a\in\M$ we have $\deg_{x_1}(a)=\deg_{x_2}(a)=3$, then $a=\sum_i \al_i\, a_i x_1^2x_2^2x_1x_2 b_i$ in $A_{3,d}$, where $\al_i\in\FF$ and $a_i,b_i\in\M_1$; 

\item[f)] $\bar{a}\bar{b}\bar{c}=-\bar{c}\bar{b}\bar{a}$ in $A_{3,d}$, where $a,b,c\in\M$. 
\end{enumerate}
\end{lemma}
\begin{proof}
{\bf a)} Let $c\in\M$. Considering $T_2(x,c)$ we obtain
\begin{eq}\label{a1}xcx=-x^2c-cx^2 \text{ in } A_{3,d},\end{eq}
and considering $T_2(x,xc)$ we obtain
\begin{eq}\label{a}xcx^2=-x^2cx \text{ in } A_{3,d}.\end{eq}
The required follows from~\Ref{a1} and~\Ref{a}.

\smallskip%
{\bf b)} See Proposition~1 of~\cite{Lopatin_Comm1}.

\smallskip%
{\bf c)} The first part follows from part~a) of Lemma~\ref{lemma_N3} and the second part is a consequence of the first part.

\smallskip%
{\bf d)} Applying part~a) of Lemma~\ref{lemma_N3}, we obtain $ab\cdot c\cdot ba=a^2 b^2c + cb^2a^2 + a^2cb^2 + b^2ca^2$ in $A_{3,d}$. To complete the proof, we consider $T_3(a^2,b^2,c)$.

\smallskip%
{\bf e)} See Statement~7 of~\cite{Lopatin_Comm1}. 

\smallskip%
{\bf f)} The equality $T(b,a,c)+T(b,c,a)-T_3(b,\bar{a},\bar{c})=0$ in $A_{3,d}$ gives the required.
\end{proof}

\begin{lemma}\label{lemma_pi}
Let $\Char{\FF}=3$ and $1\leq i\leq d$. For any homogeneous $e\in\M_{\FF}$ of multidegree $(\de_1,\ldots,\de_d)$ with $\de_i<3$ and $\de_1+\cdots+\de_{i-1}+\de_{i+1}+\cdots+\de_d>0$  we define $\pi_i(e)\in\M_{\FF}$ as the result of the substitution $x_i\to 1$, $x_i^T\to1$ in $a$, where $1$ stands for the unity of $\M_1$. 

Then $e=0$ in $A_{3,d}$ implies $\pi_i(e)=0$ in $A_{3,d}$.   
\end{lemma}
\begin{proof} Let $a,b,c\in\M_{\FF}$. By definition, $\pi_i(ab)=\pi_i(a)\pi_i(b)$. It is not difficult to see that $\pi_i(\bar{a})=\bbar{\pi_i(a)}$. 
Then by straightforward calculations we can show that 
$\pi_i(T_2(a,b))=0$, $\pi_i(T_3(a,b,c))=0$, and $\pi_i(T(a,b,c))=0$ in $A_{3,d}$. The proof is completed.  
\end{proof}

\begin{lemma}\label{lemma_four}
The equality 
$$\bar{a} u \bar{b} v \bar{c} w\bar{e}=0$$%
holds in $A_{3,d}$ for all $a,b,c,e\in\M$ and $u,v,w\in\M_1$.
\end{lemma}
\begin{proof}
In this proof all elements belong to $\M$ and all equalities are considered in $A_{3,d}$. 

Consider $T(\bar{a}\bar{b},c,e)=0$. Using  $\bbar{\bar{a}\bar{b}}=\bar{a}\bar{b}-\bar{b}\bar{a}$ and part~f) of Lemma~\ref{lemma_N3}, we obtain 
\begin{eq}\label{eq5}
\bar{a}\bar{b}\bar{c}\bar{e}=0.
\end{eq}
The equalities $T(u,a,b)\bar{c}\bar{e}=0$ and $\bar{a}\bar{b}T(u,c,e)=0$ together with~\Ref{eq5} imply
\begin{eq}\label{eq6}
\bar{a}u\bar{b}\bar{c}\bar{e}=-\bar{a}\bar{b}u\bar{c}\bar{e}=\bar{a}\bar{b}\bar{c}u\bar{e}.
\end{eq}
Hence $\bar{a}T(u,b,c)\bar{e}=0$ implies 
\begin{eq}\label{eq7}
\bar{a}u\bar{b}\bar{c}\bar{e}=\bar{a}\bar{b}u\bar{c}\bar{e}=\bar{a}\bar{b}\bar{c}u\bar{e}=0.
\end{eq}
It follows from $T_3(\bar{a}u,\bar{b}\bar{c},v\bar{e})=0$ together with~\Ref{eq5} and~\Ref{eq7} that 
\begin{eq}\label{eq8}
\bar{a}u\bar{b}\bar{c}v\bar{e}=0.
\end{eq}
Thus $\bar{a}u T_3(\bar{b},v,\bar{c}\bar{e})=0$ and $T_3(\bar{a}\bar{b},u,\bar{c})v\bar{e}=0$ imply
\begin{eq}\label{eq9}
\bar{a}u\bar{b}v\bar{c}\bar{e}=\bar{a}\bar{b}u\bar{c}v\bar{e}=0.
\end{eq}
Considering $T_3(\bar{a}u,\bar{b}v\bar{c},w\bar{e})=0$, we obtain
\begin{eq}\label{eq10}
\bar{a}u\bar{b}v\bar{c}w\bar{e}=0.
\end{eq}
The claim of the lemma is proved.
\end{proof}

\begin{lemma}\label{lemma_Akey}
Let $\Char{\FF}=3$. For $0\leq s\leq 4$ we assume that $x_1,\ldots,x_{8-2s}$ are pairwise different letters, $a_1,\ldots,a_s,v,w\in\M$, and $u_1,\ldots,u_{s+1}\in\M_1$. If $\deg_{x_i}(u_1\cdots u_{s+1})=3$ for $1\leq i\leq 8-2s$, then the following equalities hold in $A_{3,d}$: 
\begin{enumerate}
\item[a)] $u_1\bar{a}_1\cdots u_s\bar{a}_s u_{s+1}=0$;

\item[b)] a product of $vw-wv$, $\bar{a}_1,\ldots,\bar{a}_{s-1}$, $u_1,\ldots,u_{s+1}$ in any succession is equal to zero, where $s>0$.
%
%
%
%
%
%
%
%
\end{enumerate}
In particular,  if $u\in\M_1$ and $\deg_{x_i}(u)=3$ for $1\leq i\leq 8$, then $u=0$.
\end{lemma}
\begin{proof}
In this proof all equalities are considered in $A_{3,d}$. We prove by decreasing induction on $s$.

\smallskip
{\bf 1.} Let $s=4$. Then part~a) follows from Lemma~\ref{lemma_four}. 

\smallskip
{\bf 2.} Assume that part~a) holds for some $s>0$. We set $a_s=vw$ in part~a). Thus, %
$$u_1\bar{a}_1 \cdots u_{s-1}\bar{a}_{s-1} u_s (vw-wv) u_{s+1}=0.$$ %
Considering $a_i$ with $1\leq i<s$ instead of $a_s$, we complete the proof of part~b) for the given $s$.

\smallskip
{\bf 3.} Assume that part~b) holds for some $s=k+1> 0$. We claim that part~a) is valid for $s=k$. We set $W=x_1^2x_2^2x_1x_2$.  By part~e) of Lemma~\ref{lemma_N3}, $u_1\bar{a}_1\cdots u_k\bar{a}_k u_{k+1}$ is a linear combination of elements 
\begin{eq}\label{eq11}
v_1 w_1 \cdots v_{k+1} w_{k+1} v_{k+2}, 
\end{eq} %
where pairwise different $w_1,\ldots,w_{k+1}$ belong to the set $\{\bar{a}_1,\ldots,\bar{a}_k,W\}$ and $v_1,\ldots,v_{k+2}$ range over $\M_1$. For $v=x_1^2x_2^2$ and $w=x_1$ we have $W=vwx_2$ and $wvx_2=0$. Part~b) for $s=k+1$ implies that all elements from~\Ref{eq11} are zero and the claim is proved.
\end{proof}

\section{Relations for $R^{O(n)}$ of multidegree $(1,\de_2\ldots,\de_d)$}\label{section_multi}

In this section we assume that $n$ is arbitrary. 

\begin{lemma}\label{lemma1_multi}
Any homogeneous relation of $\ov{R^{O(n)}}$ of multidegree $(1,\de_2\ldots,\de_d)$ follows from the relations:
\begin{enumerate}
\item[$\bullet$] $\si_{1,t-1;r;r}(e,a;b;c)\equiv0$,

\item[$\bullet$] $\si_{t;1,r-1;r}(a;e,b;c)\equiv0$,
\end{enumerate} 
where $t+2r=n+1$, $a,b,c$ range over $\M_{\FF}$, and $e$ ranges over the elements of $\M$ such that $\deg_{x_1}(e)+\deg_{x_1^T}(e)=1$.
\end{lemma}
\begin{proof}
We assume that $f=\si_{\un{t};\un{r};\un{s}}(\un{a};\un{b};\un{c})$ is a homogeneous relation of multidegree $(1,\de_2\ldots,\de_d)$, where $|\un{t}|+2|\un{r}|>n$,  $a_i,b_j,c_k\in\M$, and consider the set of all $t_i,r_j,s_k$. Then at least one element of this set is equal to $1$. We set $u=\#\un{a}$, $v=\#\un{b}$, and $w=\#\un{c}$. 

If $t_1=1$, then applying Lemma~5.2 of~\cite{Lopatin_Orel} it is not difficult to see that the next equality holds in $\ov{\N_{\si}}$:
$$f= %
-\sum_{i=2}^u \si_{\De_i} (a_1 a_i,a_2,\ldots,a_u; \un{b}; \un{c}) %
-\sum_{j=1}^v \si_{\Theta_j} (a_2,\ldots,a_u; a_1b_j,b_1,\ldots,b_v; \un{c})+$$ %
$$+\sum_{j=1}^v \si_{\Theta_j} (a_2,\ldots,a_u; a_1b_j^T,b_1,\ldots,b_v; \un{c}),
$$
where
$\De_i=(1,t_2,\ldots,t_i -1,\ldots,t_u;\un{r};\un{s})$ and  $\Theta_j=(t_2,\ldots,t_u;1,r_1,\ldots,r_j-1,\ldots, r_v;\un{s})$.

If $r_1=1$, then the next equality holds in $\ov{\N_{\si}}$:
$$f= %
-\sum_{i=1}^u \si_{\De_i} (\un{a}; b_1a_i^T,b_2,\ldots,b_v; \un{c}) %
-\sum_{k=1}^w \si_{\Theta_k}   (b_1c_k,  a_1,\ldots,a_u; b_2,\ldots,b_v; \un{c})+$$ %
$$+\sum_{k=1}^w \si_{\Theta_k} (b_1c_k^T,a_1,\ldots,a_u; b_2,\ldots,b_v; \un{c}),
$$
where
$\De_i=(t_1,\ldots,t_i -1,\ldots,t_u;1,r_2,\ldots,r_v;\un{s})$, $\Theta_k=(1,t_1,\ldots,t_u;r_2,\ldots,r_v;s_1,\ldots,s_k -1,\ldots,s_w)$.

All other cases can be treated similarly. We repeat this procedure and use~\Ref{eq13} to obtain the required. 
\end{proof}

\section{Reduction to $A_{3,d}$}\label{section_reduction}

In this section we assume $n=3$. This section is dedicated to the proof of the following theorem and its corollary. We say that $f\in\M_{\FF}$ does not contain a letter $x$ if $f$ is a linear combination of $a_1,\ldots,a_s\in\M$, where $\deg_x a_i=0$ for all $i$.

\begin{theo}\label{theo_reduction} Any relation for $\ov{R^{O(3)}}$ follows from the relations:
\begin{enumerate}
\item[(A)] $\tr(au)\equiv0$ for $u\in\M_{\FF}$ such that $u=0$ in $A_{3,d}$;

\item[(B)] $\si_3(ab)\equiv0$;

\item[(C)] $\si_t(a)\equiv0$ for $t>3$;
\end{enumerate}
where $a,b\in\M$.
\end{theo}

\begin{lemma}\label{lemma1_reduction} 
\begin{enumerate}
\item[a)] $\si_2(a)=-\frac{1}{2}\tr(a^2)$ and $3\si_3(a)=\tr(a^3)$ in $\ov{\N_{\si}}$, where $a\in\M_{\FF}$.

\item[b)] Any relation $h\equiv0$ of multidegree $(1,\de_2,\ldots,\de_d)$ follows from (A). 

\item[c)] Any relation $\si_{t,r}(a,b,c)\equiv0$, where $t+2r>3$, $a,b,c\in\M_{\FF}$, and $t$ or $r$ is not divisible by $\Char{\FF}$ (or $\Char{\FF}=0$), follows from (A). 

\item[d)] If $a\in\M$ and $\deg_{x_1}{a}>3$, then $\tr(a)\equiv0$ follows from (A).

\item[e)] If $\Char{\FF}=3$ and $u=0$ in $A_{3,d}$, then $\tr(u)\equiv0$ follows from  (A).
\end{enumerate}
\end{lemma}
\begin{proof} 
{\bf a)} Since $\tr(a^2)=\tr(a)^2-2\si_2(a)$ in $\N_{\si}$ (see Lemma~3.1 from~\cite{Lopatin_Orel}), we obtain the first equality. The second equality can be proved similarly.

{\bf b)} 
By Lemma~\ref{lemma1_multi}, it is enough to show that $f_1=\si_{1,t-1;r;r}(e,a;b;c)\equiv0$ and $f_2=\si_{t;1,r-1;r}(a;e,b;c)\equiv0$ with $t+2r=4$ follow from (A), where $a,b,c\in \M_{\FF}$ and $e\in \M$. To complete the proof we consider the following cases. 
 
If $t=2$ and $r=1$, then $f_1=-\tr(e T(a,b,c))$ and $f_2= -\tr(e T_2(a,\bar{c})) + \tr(e T(a,a,c))$ in $\ov{\N_{\si}}$.

If $t=4$ and $r=0$, then $f_1=-\tr(e a^3)$ in $\ov{\N_{\si}}$.

If $t=0$ and $r=2$, then $f_2=-\tr(e \bar{c} \bar{b} \bar{c})=-\tr(e T(\bar{c},b,c))$ in $\ov{\N_{\si}}$.

\smallskip
{\bf c)} If $t$ is not divisible by $\Char{\FF}$, then $\si_{t,r}(a,b,c)=\frac{1}{t}\si_{t-1,1;r;r}(a,a;b;c)$ (see Lemma~5.5 of~\cite{Lopatin_Orel}). Part~b) concludes the proof. The other case is alike. 

\smallskip%
{\bf d)} It follows from the proof of part~a) of Lemma~\ref{lemma_N3} and the linearity of $\tr$.

\smallskip%
{\bf e)} For $a,b,c\in\M$ we have $\tr(T_2(a,b))=3\tr(a^2b)$, $\tr(T_3(a,b,c))=3\tr(abc)+3\tr(acb)$, and $\tr(T(a,b,c))=3\tr(a\bar{b}\bar{c})$ in $\N_{\si}$. By part~a) of Lemma~\ref{lemma1_reduction}, we have $\tr(T_1(a))=3\si_3(a)$ in $\ov{\N_{\si}}$. The proof is completed. 
\end{proof}

\begin{lemma}\label{lemma_reduction_verification}
We have that (A) and (B) are relations for $\ov{R^{O(3)}}$. 
\end{lemma}
\begin{proof}
The proof of part~b) of Lemma~\ref{lemma1_reduction} implies that $\tr(e a^3)\equiv0$ and $\tr(e T(a,b,c))\equiv0$ in $\ov{R^{O(3)}}$ for all $a,b,c\in\M_{\FF}$ and $e\in\M$. The fact that  $\tr$ is linear shows that (A) is a relation. 

Since $\si_3(AB)=\si_3(A)\si_3(B)$ for all $3\times 3$ matrices $A$ and $B$ over the polynomial  ring $R$ (see Section~\ref{section_intro}), we obtain $\si_3(ab)\equiv0$.
\end{proof}

\begin{lemma}\label{lemma2_reduction} The relations $\si_{0,3}(v,w)\equiv0$ and $\si_{3,3}(u,v,w)\equiv0$, where $u,v,w\in\M_{\FF}$, follow from (A), (B), and (C).
\end{lemma}
\begin{proof}
Clearly, without loss of generality we can assume that $u=x_1$, $v=x_2$, $w=x_3$ are letters and $d\geq3$.  Consider $t,r\geq0$ such that $t+2r>3$ and $t$ or $r$ is odd.

Let $x\in\M$ be a letter and $a=b_0a_1b_1\cdots a_p b_p \in\M$, where $a_1,\ldots,a_p\in\M$ are words in $x,x^T$ and $b_1,\ldots,b_{p-1}\in\M$, $b_0,b_p\in\M_1$ do not contain $x,x^T$, i.e.,  $\deg_{x}(b_i)=\deg_{x^T}(b_i)=0$ for $0\leq i\leq p$. We say that $a$ is {\it fixed} by $x$ if there is no non-trivial cyclic permutation $\pi\in S_p$ such that  
$$a_1 = a_{\pi(1)},\ldots, a_p = a_{\pi(p)}\,\text{ or }\, 
a_1^T = a_{\pi(p)},\ldots, a_p^T = a_{\pi(1)}.
$$
Assume that $a=b_0a_1b_1\cdots a_p b_p$ is fixed by $x$ and $b_0\cdots b_p=y_1\cdots y_s$, where $y_i$ is a letter for all $i$. Then the following elements 
$$c_0a_1c_1\cdots a_p c_p\, ,$$ %
where $c_1,\ldots,c_{p-1}\in\M$, $c_0,c_p\in\M_1$, and $c_0\cdots c_p=z_{\rho(1)}\cdots z_{\rho(s)}$ for $z_i\in\{y_i,y_i^T\}$ ($1\leq i\leq s$) and $\rho\in S_s$, are primitive and pairwise different with respect to the $\sim$-equivalence.

For $x\in\{x_1,x_2,x_3\}$ we denote $Y=\{x_2,x_3\}\setminus\{x\}$. The above mentioned property together with the definition of $\si_{t,r}$ implies that  
\begin{eq}\label{eq1}
\si_{t,r}(x_1,x_2,x_3)=f_1+f_2+f_3 \text{ in }\ov{\N_{\si}} .
\end{eq}
Here 
\begin{enumerate}
\item[$\bullet$] $f_1=\sum (-1)^{\xi_{a}}\tr(a)|_{y\to \bar{y},\;y\in Y}$, where the sum ranges over $a\in\ov{I_{t,r}}$ such that $a$ is fixed by $x$ 
and $\deg_{y}(a)=r$ for all $y\in Y$;

\item[$\bullet$] $f_2=\sum (-1)^{\xi_{a}}\tr(a)$, where the sum ranges over $a\in\ov{I_{t,r}}$ such that $a$ is not fixed by $x$;

\item[$\bullet$] $f_3\equiv0$ follows from (B) and (C).
\end{enumerate} 
Let us consider $\si_{t,r}$ from the formulation of the lemma.

\smallskip %
{\bf a)} Let $t=0$ and $r=3$. We take $x=x_2$. By part~c) of Lemma~\ref{lemma_N3},
$$f_1=\tr(x_2 \bar{x}_3 x_2 \bar{x}_3 x_2^T \bar{x}_3)\equiv \tr(\bar{x}_3^2 x_2^2 x_2^T \bar{x}_3)\equiv0 \text{ and }$$
$$f_2= \tr(x_2 x_3 x_2 x_3 x_2 x_3^T) - \tr(x_2 x_3 x_2 x_3^T x_2 x_3^T)\equiv %
\tr(x_3^2 x_2^2 x_2 x_3^T) - \tr(x_2 x_3 (x_3^T)^2 x_2^2) \equiv0
$$
follow from (A).

\smallskip %
{\bf b)} Let $t=r=3$. We take $x=x_1$. By Lemma~\ref{lemma_four}, $f_1\equiv0$ follows from (A). Clearly, $a\in\ov{I_{t,r}}$ is not fixed by $x$ if and only if $a\sim x a_1 x a_2 x a_3$, where $a_1,a_2,a_3\in\M$.  Thus, $f_2=h_1+h_2$ in $\ov{\N_{\si}}$, where
$$
\begin{array}{rcl}
h_1=\tr(x_1x_2^Tx_3^T x_1x_2^Tx_3 x_1x_2x_3^T)&+ %
&\tr(x_1x_2^Tx_3^T x_1x_2x_3^T x_1x_2^Tx_3)\\ 
-\tr(x_1x_2^Tx_3^T x_1x_2^Tx_3 x_1x_2x_3) & - &\tr(x_1x_2^Tx_3^T x_1x_2x_3^T  x_1x_2x_3) \\
-\tr(x_1x_2^Tx_3^T x_1x_2x_3   x_1x_2^Tx_3)&-&\tr(x_1x_2^Tx_3^T x_1x_2x_3   x_1x_2x_3^T)\\ 
+\tr(x_1x_2^Tx_3   x_1x_2x_3^T x_1x_2x_3)&  +&\tr(x_1x_2^Tx_3   x_1x_2x_3   x_1x_2x_3^T)\\
\end{array} 
$$ %
and $h_2$ is a linear combination of elements $h=\tr(x_1 y_2 y_3 x_1 y_2 y_3 x_1 z_2 z_3)$, where $y_i,z_i\in\{x_i,x_i^T\}$ for $i=2,3$. Since $x_1 y_2 y_3\cdot x_1 y_2 y_3\cdot x_1 =  y_3^2 y_2^2 x_1^2\cdot x_1=0$ in $A_{3,d}$ (see part~c) of Lemma~\ref{lemma_N3}), $h\equiv0$ follows from (A). 

Using the equality $b^T=b-\bar{b}$, we can get rid of $x_2^T$ and $x_3^T$ in $h_1$. Hence 
$$
\begin{array}{rcl}
h_1=\tr(x_1\bar{x}_2\bar{x}_3 x_1\bar{x}_2x_3 x_1x_2\bar{x}_3)&+    %
&\tr(x_1\bar{x}_2\bar{x}_3 x_1x_2\bar{x}_3 x_1\bar{x}_2x_3)\\ 
+\tr(x_1\bar{x}_2x_3   x_1x_2\bar{x}_3 x_1x_2x_3)&  +&\tr(x_1\bar{x}_2x_3   x_1x_2x_3   x_1x_2\bar{x}_3)+h_3\\
\end{array} 
$$ %
in $\ov{\N_{\si}}$, where $h_3$ is a linear combination of elements $h$ as above. By Lemma~\ref{lemma_four}, 
$$\tr(x_1\bar{x}_2\bar{x}_3 x_1\bar{x}_2x_3 x_1x_2\bar{x}_3)+     
\tr(x_1\bar{x}_2\bar{x}_3 x_1x_2\bar{x}_3 x_1\bar{x}_2x_3)\equiv0$$ %
follows from (A). Applying part~d) of Lemma~\ref{lemma_N3} to 
$x_3\cdot x_1x_2\cdot\bar{x}_3\cdot x_1x_2 \cdot x_3$ and using part~c) of Lemma~\ref{lemma_N3}, we obtain 
$$\tr(x_1\bar{x}_2x_3   x_1x_2\bar{x}_3 x_1x_2x_3) + \tr(x_1\bar{x}_2x_3   x_1x_2x_3   x_1x_2\bar{x}_3)\equiv -\tr(x_1 \bar{x}_2 \bar{x}_3 x_3^2 x_2^2 x_1^2)\equiv0$$
follows from (A). The proof is completed.
\end{proof}

%
%
%
%

\begin{proof_without_dot} {\it $\!\!\!\!$ of Theorem~\ref{theo_reduction}.} 
By Theorem~\ref{theo_rel_main} and Lemma~\ref{lemma_reduction_verification}, it is enough to show that $f=\si_{t,r}(a,b,c)\equiv0$, where $t+2r>3$ and $a,b,c\in\M_{\FF}$, follows from (A), (B), and (C).

Applying part~a) of Lemma~\ref{lemma1_reduction}, we obtain that $f=f_1+f_2+f_3$ in $\ov{\N_{\si}}$, where  $f_1=\sum_i\al_i\tr(a_i)$, $f_2=\sum_j\be_j\si_3(b_j)$, $\al_i,\be_j\in \FF$, $a_i,b_j\in\M$, and $f_3$ follows from (C). Since $\deg{f}>3$, $f_2\equiv0$ follows from (B).

If $t>6$ or $r>6$, then $f_1\equiv0$ follows from (A). To prove this claim, note that if $t>6$, then $\deg_{a}{a_i}>3$ or $\deg_{a^T}{a_i}>3$, where we consider $a_i$ as a word in $a,b,c,a^T,b^T,c^T$. Part~d) of Lemma~\ref{lemma1_reduction} implies that $\tr(a_i)\equiv0$ follows from (A).


\smallskip
{\bf a)} Let $\Char{\FF}=0$ or $\Char{\FF}\geq5$. Part~b) of Lemma~\ref{lemma_N3} implies that if $\deg{f}>6$, then $f_1\equiv0$ follows from (A). If $\deg{f}\leq6$, then $f\equiv0$ follows from (A) by part~c) of Lemma~\ref{lemma1_reduction} or $f=\si_5(x)\equiv0$ follows from (C), where $x$ is a letter. 

\smallskip
{\bf b)} Let $\Char{\FF}=3$. In the rest of the proof we apply part~e) of Lemma~\ref{lemma1_reduction} without reference to it. We have $f_1=\sum_k\ga_k \tr(c_k)$, where $\ga_k\in\FF$ and $c_k$ is a monomial in $a,b,c,\bar{a},\bar{b},\bar{c}$ for all $k$. 

We assume $r=6$. If $\deg_{b}{c_k}=\deg_{c}{c_k}=\deg_{\bar{b}}{c_k}=\deg_{\bar{c}}{c_k}=3$, then $\tr(c_k)\equiv0$ follows from (A) by Lemma~\ref{lemma_four}; otherwise, $\tr(c_k)\equiv0$ follows from (A) by part~d) of Lemma~\ref{lemma1_reduction}.

We claim that if $t=6$ and $r=3$, then $\tr(c_k)\equiv0$ follows from (A) for all $k$. Assume that this claim does not hold. Then $\deg_a(c_k)=\deg_{\bar{a}}(c_k)=3$ by part~d) of  Lemma~\ref{lemma1_reduction} and $\deg_b(c_k)=\deg_c(c_k)=3$ by Lemma~\ref{lemma_four}. Part~e) of Lemma~\ref{lemma_N3} implies that %
$c_k=\al\, W_{a\bar{a}} W_{bc}+\be\, W_{bc} W_{a\bar{a}}$ in $A_{3,d}$, where $\al,\be\in\FF$ and $W_{uv}=u^2v^2uv$ for any $u,v\in\M_{\FF}$. 
We set $v=b^2c^2$, $w=b$, and $s=4$ in part~b) of Lemma~\ref{lemma_Akey} and obtain  $W_{a\bar{a}} W_{bc}= W_{bc} W_{a\bar{a}}=0$ in $A_{3,d}$. Hence $\tr(c_k)\equiv0$ follows from (A); a contradiction. Thus the claim is proved.

If $t=6$ and $r=0$, then it is not difficult to see that $f=\si_6(a)\equiv0$ follows from (A), (B), and (C) by Lemma~9 of~\cite{Lopatin_Comm1}. 

Part~c) of Lemma~\ref{lemma1_reduction} and Lemma~\ref{lemma2_reduction} complete the consideration of the case of $\Char{\FF}=3$. 
\end{proof_without_dot}

\begin{cor}\label{cor1}
Let $u,v\in\M_{\FF}$ do not contain $x_1,x_1^T$. We have 
\begin{enumerate}
\item[a)] $\tr(ux_1)\equiv0$ if and only if $u=0$ in $A_{3,d}$;

\item[b)] if $\Char{\FF}=3$, then $\tr(ux_1^2)\equiv0$ if and only if $u=0$ in $A_{3,d}$.
%
\end{enumerate}
\end{cor}
\begin{proof}
{\bf a)}  If $u=0$ in $A_{3,d}$, then $\tr(ux_1)\equiv0$ by Lemma~\ref{lemma_reduction_verification}.

If $\tr(ux_1)\equiv0$, then $\tr(ux_1)=\sum_i\al_i \tr(u_ia_i)$ in $\ov{\N_{\si}}$, where $\al_i\in\FF$, $a_i\in\M$, $u_i\in\M_{\FF}$ is $\NN^d$-homogeneous with $u_i=0$ in $A_{3,d}$, and $\deg_{x_1}(a_iu_i)+\deg_{x_1^T}(a_iu_i)=1$ (see Theorem~\ref{theo_reduction}). For $a,b,c,e\in\M_{\FF}$ the following equalities in $\N_{\si}$ hold:
$$\tr(T_2(a,b)\,e)=\tr(T_2(a,e)\,b),$$
$$\tr(T_3(a,b,c)\,e)=\tr(T_3(e,b,c)\,a),$$
$$\tr(T(a,b,c)\,e)=\tr(T(e,b,c)\,a),$$
$$\tr(T(a,b,c)\,e)=\tr\left( (T_3(a,e,\bar{c}) - T(a,e,c) - T(e,a,c)) \,b \right),$$
$$\tr(T(a,b,c)\,e)=\tr\left( (T_3(a,\bar{b},e) - T(a,b,e) - T(e,b,a)) \,c \right).$$
Thus $\tr(ux_1)=\sum_j\be_j \tr(v_jx_1)$ in $\ov{\N_{\si}}$, where $\be_j\in\FF$ and $v_j\in\M_{\FF}$ satisfies $v_j=0$ in $A_{3,d}$. Therefore $u=\sum_j\be_j v_j$ in $\M_{\FF}$ and the proof is completed. 

\smallskip
{\bf b)} If $\tr(ux_1^2)\equiv0$, then $\tr(u\,(x_1+x_{d+1})^2)\equiv0$, where $x_{d+1}$ stands for a new letter. Taking the homogeneous components of multidegrees $(1,\de_2,\ldots,\de_d,1)$ we obtain $\tr((a x_1+x_1 a) x_{d+1})\equiv0$. By part~a), $u x_1 + x_1 u=0$ in $A_{3,d}$. The application of $\pi_1$ from Lemma~\ref{lemma_pi} completes the proof.  
%
\end{proof}

\section{$D_{max}$ and the nilpotency degree of $A_{3,d}$}\label{section_proof}

Theorem~\ref{theo_appl} is a consequence of Lemma~\ref{lemma3_proof} which is proved in this section. For an $n\times n$ matrix $X$ we denote $\bar{X}=X-X^T$.

\begin{lemma}\label{lemma1_proof}
The invariant $\tr(X_1^2\bar{X}_1^2X_1\bar{X}_1)\in R^{O(3)}$ is indecomposable, i.e., $\tr(x_1^2\bar{x}_1^2x_1\bar{x}_1)\not\equiv0$.
\end{lemma}
\begin{proof} For short, we write $X$ instead of $X_1$. Assume that the claim of this lemma  does not hold. By part~d) of Lemma~\ref{lemma1_reduction}, any indecomposable element of $R^{O(3)}$ of multidegree $(k)$, where $k<6$, is a linear combination of the following elements:
\begin{eq}\label{eq3}
\tr(X),\; \si_2(X),\; \si_3(X),\; \tr(X X^T),\; \tr(X^2 X^T),\; \tr(X^2(X^T)^2),\; \si_2(X X^T).
\end{eq}
Then 
\begin{eq}\label{eq4}
\tr(X^2\bar{X}^2 X\bar{X})=\sum_i \al_i f_{i_1}\cdots f_{i_q}, 
\end{eq} %
where $\al_i\in\FF$ and $f_j$ is an element from~\Ref{eq3}. Moreover, if we take an arbitrary $3\times 3$ matrix over $R$ instead of $X$, then~\Ref{eq3} remains valid and the coefficients $\al_i$ do not depend on $X$. Denote by $\al,\be,\ga\in \FF$, respectively, the coefficients of $\tr(XX^T)^3$, $\tr(XX^T)\tr(X^2 (X^T)^2)$, and $\tr(XX^T) \si_2(XX^T)$, respectively, in~\Ref{eq4}. For $a,b\in R$ we set %
$$X=\left(
\begin{array}{ccc}
0 & a& 0\\
0 & 0& b\\
0 & 0& 0\\
\end{array}
\right). %
$$
If we take $b=0$, then $\al=0$. If we assume $a,b\neq0$, then 
$$\al\tr(XX^T)^3 + \be\tr(XX^T)\tr(X^2 (X^T)^2)+ \ga\tr(XX^T) \si_2(XX^T)=0.$$
Thus
$$a^4b^2 (1+\be+\ga) + a^2b^4(\be+\ga)=0$$
for all $a,b\in R$; a contradiction.
\end{proof}

\begin{lemma}\label{lemma_DA} Let $\Char{\FF}=3$. Then
\begin{enumerate}
\item[a)] for any $a\in\M$ satisfying $a\neq0$ in $A_{3,d}$ we have $\deg{a}\leq 2d+7$.

\item[b)] there exists an $a\in\M$ such that $a\neq0$ in $A_{3,d}$ and $\deg{a}= 2d+4$.
\end{enumerate}
\end{lemma}
\begin{proof}
{\bf a)} In this proof all equalities are considered in $A_{3,d}$. We can assume that $a$ is a word in $x_1,\ldots,x_d$, $\bar{x}_1,\ldots,\bar{x}_d$. By part~a) of Lemma~\ref{lemma_N3}, $\deg_{x_i}(a)\leq3$ for all $i$. Let $k=\deg_{\ov{x}_1}(a)+\cdots+\deg_{\ov{x}_d}(a)$. By Lemma~\ref{lemma_four}, $k<4$.
Moreover, by part~a) of Lemma~\ref{lemma_Akey} there are no pairwise different $1\leq i_1,\ldots,i_{8-2k}\leq d$ such that $\deg_{x_{i_j}}(a)=3$ for all $1\leq j\leq 8-2k$. Thus, $\deg{a}\leq k + 2d + (7-2k)\leq 2d+7$.

\smallskip
{\bf b)} Consider $a_d=x_1^2 \bar{x}_1^2 x_1 \bar{x}_1 x_2^2 \cdots x_d^2\in\M$. If $d=1$, then $a_d\neq0$ in $A_{3,d}$ by part~e) of Lemma~\ref{lemma1_reduction} and Lemma~\ref{lemma1_proof}. We assume that $d>1$. If $a_d=0$ in $A_{3,d}$, then applying $\pi_d,\ldots,\pi_2$ from Lemma~\ref{lemma_pi} to $a_d$ we obtain that $a_{1}=0$ in $A_{3,d}$; a contradiction. 
\end{proof}

\begin{lemma}\label{lemma3_proof}
\begin{enumerate}
\item[a)] If $\Char{\FF}=3$ and $d\geq1$, then $2d+4\leq D_{\rm max}\leq 2d+7$. 

\item[b)] If $\Char{\FF}\neq 2$ and $d=1$, then $D_{\rm max}=6$.

\item[c)] If $\Char{\FF}\neq 2,3$ and $d\geq1$, then $D_{\rm max}=6$.
\end{enumerate}
\end{lemma}
\begin{proof}
{\bf a)} Let $\tr(a)\not\equiv0$, where $a\in\M_{\FF}$. By part~e) of Lemma~\ref{lemma1_reduction}, $a\neq0$ in $A_{3,d}$. Part~a) of Lemma~\ref{lemma_DA} implies $\deg{a}\leq 2d+7$. Thus part~a) of Lemma~\ref{lemma1_reduction} and equality~(B) from Theorem~\ref{theo_reduction} show that $D_{\rm max}\leq 2d+7$.

Consider $a_d=x_1^2 \bar{x}_1^2 x_1 \bar{x}_1 x_2^2 \cdots x_d^2\in\M$. If $d=1$, then $\tr(a_d)\not\equiv0$ by Lemma~\ref{lemma1_proof}. We assume that $d>1$. If $\tr(a_d)=0$, then $a_{d-1}=0$ in $A_{3,d}$ by part~b) of Corollary~\ref{cor1}; a contradiction with the proof of part~b) of Lemma~\ref{lemma_DA}.  

\smallskip
{\bf b)} Let $a\in\M_{\FF}$ and $d=1$. If $\deg{a}>6$, then $\deg_{x_1}(a)>3$ or $\deg_{x_1^T}(a)>3$; and part~a) of Lemma~\ref{lemma_N3} implies $a=0$ in $A_{3,d}$. Lemma~\ref{lemma1_proof} completes the proof.

\smallskip
{\bf c)} The required follows from part~b) of Lemma~\ref{lemma_N3}, equality~(A) of Lemma~\ref{theo_reduction}, and Lemma~\ref{lemma1_proof}.
\end{proof}

\remark{} Denote by $D_{\rm nil}$ the {\it nilpotency degree} of $A_{3,d}$, i.e., the minimal $s$ such that $a_1\cdots a_s=0$ in $A_{3,d}$ for all $a_1,\ldots,a_s\in \M_{\FF}$. It is not difficult to see that if $d>1$, then the same estimations as for $D_{\rm max}$ in Lemma~\ref{lemma3_proof} are also valid for $D_{\rm nil}-1$.

\conj\label{conj1} If $\Char{\FF}=3$ and $d\geq7$, then $D_{\rm max}=2d+7$.

\remark\label{} Let $\Char{\FF}=3$ and $d\geq7$. To prove Conjecture~\ref{conj1} it is enough to show that 
\begin{eq}\label{eq12}
x_1^2 W_{23} x_1 W_{45} W_{67}\neq0 \text{ in }A_{3,d},
\end{eq}
where $W_{ij}=x_i^2x_j^2x_ix_j$. This claim follows from part~b) of Corollary~\ref{cor1} and Statement~8 of~\cite{Lopatin_Comm1}. 


\end{document}